\documentclass{amsart}
\usepackage{amssymb}
\usepackage{latexsym}
\usepackage{amsmath}
\usepackage{euscript}
\usepackage{graphics}
\usepackage[all]{xy}
\usepackage{amsfonts,amsthm}
\usepackage[centering]{geometry}
\usepackage{MnSymbol} 
\usepackage{color}

\newcommand{\be}{\begin{equation}}
\newcommand{\ee}{\end{equation}}
\newcommand{\ba}{\begin{eqnarray}}
\newcommand{\ea}{\end{eqnarray}}
\newcommand{\baa}{\begin{eqnarray*}}
\newcommand{\eaa}{\end{eqnarray*}}
\newcommand{\bb}{\begin{thebibliography}}
\newcommand{\eb}{\end{thebibliography}}

\newcommand{\bi}[1]{\bibitem{#1}}
\newcommand{\lab}[1]{\label{#1}}
\newcommand{\re}[1]{(\ref{#1})}

\newcounter{my}
\newcommand{\he}%
   {\stepcounter{equation}\setcounter{my}%
   {\value{equation}}\setcounter{equation}0%
   }%
\newcommand{\she}%
   {\setcounter{equation}{\value{my}}%
    }%

\renewcommand\t{\tilde}

\newcommand\ve{\varepsilon}

\newtheorem{theorem}{Theorem}[section]
\newtheorem{pr}[theorem]{Proposition}

\newtheorem{lemma}[theorem]{Lemma}

\theoremstyle{definition}

\numberwithin{equation}{section}

\begin{document}


\title[Persymmetric OPUC]{Mirror symmetric polynomials orthogonal on the unit circle}


\author{Alexei Zhedanov}

\address{School of Mathematics, Renmin University of China, Beijing 100872,
China}

\begin{abstract}
We introduce and study a special family of polynomials orthogonal on the unit circle (OPUC). These OPUC satisfy a mirror symmetry property of their Verblunsky coefficients. Several equivalent conditions for the OPUC to be mirror symmetric are presented. Corresponding unitary CMV matrices satisfy simple algebraic relations similar to relations  for persymmetric tridiagonal matrices. We present three explicit examples of mirror symmetric OPUC.   
\end{abstract}



\maketitle

\section{Introduction}
\setcounter{equation}{0}
Tridiagonal (Jacobi) matrices is a common subject of research in many branches of mathematics and mathematical physics.  The reason of this is the well known Stones's theorem \cite{AG} that a selfadjoint operator in a Hilbert space with simple spectrum  can be represented as a symmetric semi-infinite Jacobi operator in a suitable orthonormal basis. Hence the Jacobi matrices provide a convenient model for studying properties of self-adjoint operators.

Another reason is that the eigenvalue problem for Jacobi matrices is related to the theory of orthogonal polynomials on the rel line (OPRL) \cite{Ismail}. The latter appear in many problems of applied mathematics and theoretical physics. Moreover, there is a family of the "classical" orthogonal polynomials \cite{Ismail}. Many  formulas for these polynomials become explicit which leads to numerous explicit examples of Jacobi matrices with prescribed properties.  

There is a special class of finite Jacobi matrices which is called the {\it persymmetric}, or mirror-symmetric, matrices \cite{Gladwell}, \cite{per}. These Jacobi matrices are "double" symmetric, i.e. they are symmetric not only with respect to the main diagonal but also with respect to the main antidiagonal. Many properties of persymmetric matrices make them a very convenient tool in applications. One could mention mechanical oscillations  \cite{Gladwell} and the problem of perfect state transfer in quantum information \cite{Albanese}, \cite{VZ_PST}.

The inverse spectral problem for Jacobi matrices is also very important in numerous applications \cite{Gladwell}. This problem can be formulated as follows. Assume that the spectrum (i.e. a set of distinct real numbers) of a finite Jacobi matrix $J$ is given. How can we reconstruct the matrix $J$?   The answer is that for generic case the spectrum is not enough to reconstruct the Jacobi matrix $J$ uniquely. We need also additional spectral data. These data can be chosen in different ways which leads to different types of the inverse spectral problem \cite{BG}.

The persymmetric Jacobi matrices provide a nice exception: these matrices are reconstructed {\it uniquely} if their spectrum is given \cite{BG}, \cite{Gladwell}. This property is very useful e.g. in construction of a model of the quantum register \cite{VZ_PST}.

Another class of operators which is common in many applications (e.g. in quantum theory) is the class of the unitary operators. There is an analog of Ston'e theorem for the unitary operators. Namely, one can show that such operators can be represented by a special class of "zigzag" pentadiagonal matrices. These matrices are called the CMV matrices \cite{CMV}, \cite{Simon}. The inverse spectral problem for CMV matrices was considered in \cite{GoKu}, \cite{GoKu2} where several useful algorithms to reconstruct CMV matrices from spectral data were presented.

Remarkably the eigenvalue problem for the CMV matrices leads to the theory of polynomials orthogonal on the unit circle (OPUC)  introduced by Szeg\H{o} around 100 years ago \cite{Szego}. In the case of finite CMV matrices one can introduce the mirror dual system of OPUC which was first considered by Delsarte and Genin in \cite{DG1}, \cite{DG2}. This duality leads to an interesting new interpretation of the famous Poncelet problem in projective geometry (see \cite{Mar_Si} for details).

The main subject of the present paper is further analysis of the mirror duality property of OPUC and of corresponding CMV matrices. We introduce the persymmetric OPUC and derive their characterization properties. It appears that many basic properties are similar to those of persymmetric OPRL. Moreover, we show how mirror duality property can be presented in terms of CMV matrices. Similarity with the case of the Jacobi persymmetric matrices takes place as well. However there are several points where the theory of persymmetric CMV matrices looks quite different. 

The paper is organized as follows. In Section 2, we recall basic properties of persymmetric polynomials on the real line (OPRL). In Section 3 and 4, general properties of the polynomials orthogonal on the unit circle (OPUC) and corresponding CMV matrices are briefly recalled. In Section 5, definitions of mirror dual and persymmetric OPUC are  presented together with several necessary and sufficient conditions for OPUC to be persymmetric. In Section 6, we present several explicit examples of persymmetric OPUC: "free" ones, persymmetric analogs of the single-moment OPUC and the Krawtchouk-type OPUC. The last example is new and nontrivial because it contains an arbitrary complex parameter which is absent for the case of the ordinary Krawtchouk polynomials. Section 7 deals with invariant properties of persymmetic CMV matrices with respect to (quasi) reflection operators. Finally, in Conclusion, we propose open problems for further research.

\section{Persymmetric Jacobi matrices and polynomials orthogonal on the real line}
\setcounter{equation}{0}
Consider the tri-diagonal or Jacobi Hermitian matrix of size $(N+1) \times (N+1)$
\be
J =
 \begin{pmatrix}
  b_{0} & g_1 & 0 &    \\
  g_{1} & b_{1} & g_2 & 0  \\
   0  &  g_2 & b_2 & g_3    \\
   &   &  \ddots &    \ddots  \\
      & & \dots &g_{N-1} & b_{N-1} & g_N \\   
            & & \dots &0 & g_N & b_N
         
\end{pmatrix}.
\lab{J_Jac} \ee 
In what follows we will assume that all entries $g_i, b_i, \: i=0, 1, \dots$ are real and moreover that $g_n > 0, \: n=1,2 ,\dots, N$. 

With the matrix $J$ one can associate a family of finite monic orthogonal polynomials defined by the recurrence relation 
\be
P_{n+1}(x) + b_n P_n(x) + u_n P_{n-1}(x) = x P_n(x), \; P_{-1}=0, P_0=1,  \lab{3_term_P} \ee
where $u_n = g_n^2>0$.
Under the conditions $u_n>0$ (which is automatic in our case) the roots of all orthogonal polynomials $P_n(x)$ are real and simple. The "final" polynomial
\be
P_{N+1}(x) = (x-x_0) (x-x_1) \dots (x-x_N), \lab{P_N+1_roots} \ee
has the real and distinct roots $x_0, x_1, \dots, x_N$. These roots coincide with the eigenvalues of the matrix $J$ \cite{Chi}. Equivalently, the characteristic polynomial of the Jacobi matrix \re{J_Jac} coincides with the polynomial $P_{N+1}(x)$.

Then the following orthogonality relation \cite{Chi}
\be
\sum_{s=0}^N w_s P_n(x_s) P_m(x_s) = h_n \: \delta_{nm}.  \lab{ort_h} \ee  
holds, where
\be
w_s = \frac{h_N}{P_{N}(x_s) P_{N+1}'(x_s)}, \quad s=0,1,\dots,N,  \lab{w_s_expl} \ee
are positive weights $w_s>0, s=0,1,\dots, N$ and where 
\be
h_n = u_1 u_2 \dots u_n \lab{h_def} \ee
are (positive) normalization constants. 
Notice that the weights are normalized
\be
\sum_{s=0}^N w_s =1. \lab{norm_w} \ee
We thus deal with polynomials orthogonal on a finite set of points of the real line (OPRL for brevity).

In what follows we will assume that the eigenvalues are 
ordered by increasing:
\be
x_0<x_1<x_2<\dots <x_N \lab{order_x}. \ee

There is an important interlacing property of the zeros of OPRL
 \cite{Atkinson}: all zeros of the polynomial $P_n(x), \: n=1,2,\dots, N$ are distinct and moreover any zero of the polynomial $P_n(x)$ lies between two neighbor zeros of the polynomial $P_{n+1}(x)$. In particular, this means that the sequence $P_N(x_s), \: s=0,1,2,\dots$ has alternating signs:
\be
P_N(x_s) = (-1)^{N+s} \left|  P_N(x_s) \right|, \quad s=0,1,\dots, N, \lab{alt_P_N}   \ee

Let $R$ be the ``reflection'' matrix
\[
R=\begin{pmatrix}
  0 & 0 & \dots & 0 & 1    \\
  0 & 0 & \dots  & 1 & 0  \\
    \dots  & \dots & \dots & \dots & \dots      \\
   1 & 0 &  \dots & 0 &0  \\
\end{pmatrix}.
\]
Clearly, $R$ is an involution, i.e.
\be
R^2 = I, \lab{R^2} \ee
where $I$ is the identity matrix. Hence the eigenvalues of the matrix $R$ are either 1 or -1.

Introduce the Jacobi matrix
\be
\hat J = R J R .
\lab{hat_J} \ee
The matrix $\hat J$  is called the {\it mirror dual} with respect to $J$.

It is easy to see that the entries $\hat g_n, \hat b_n$ of the mirror dual matrix $\hat J$ are
\be
\hat g_n = g_{N+1-n}, \: \hat b_n = b_{N-n}, \; n=0,1,\dots, N.
\lab{hat_gb} \ee
We denote $\hat P_n(x)$ the corresponding {\it mirror dual} OPRL. More exactly, the polynomials $\hat P_n(x)$ are defined by recurrence relation \re{3_term_P} where the recurrence coefficients $g_n, b_n$ are replaced with $\hat g_n, \hat b_n$.

Relation \re{hat_J} means that the matrix $\hat J$ is obtained from the matrix $J$ by a similarity transformation. Hence the spectrum $x_s, \, s=0,1,\dots, N$ of the matrix $\hat J$ is the same as of the matrix $J$.

We can thus  write down the orthogonality relation for the mirror dual polynomials
\be
\sum_{s=0}^N \hat w_s \hat P_n(x_s) \hat P_m(x_s) = \hat h_n \: \delta_{nm},  \lab{ort_P_hat} \ee  
which is similar to \re{ort_h}.

There is the important relation between the weights \cite{BG}, \cite{BS}, \cite{Bor}, \cite{VZ1}
\be
\hat w_s w_s = \frac{h_N}{{P'}^2_{N+1}(x)}.
\lab{ww_P} \ee

The symmetric Jacobi matrix $J$ is called the {\it persymmetric} (or mirror-symmetric) \cite{Gladwell} if the mirror dual matrix coincides with initial matrix: $\hat J =J$. This means that $J$ commutes with the reflection matrix
\be
JR = RJ \lab{JR}. \ee
Condition \re{JR} is equivalent to the following relations 
\be
g_{N+1-i}=g_i, \quad b_{N-i} = b_i, \quad i=0,1,\dots N. \lab{sym_ab} \ee
Hence the persymmetric matrix $J$ takes the form
\[
J =
 \begin{pmatrix}
  b_{0} & a_1 & 0 &    \\
  a_{1} & b_{1} & a_2 & 0  \\
   0  &  a_2 & b_2 & a_3    \\
   &   &  \ddots &    \ddots  \\
      & & \dots &a_2 & b_1 & a_1 \\   
            & & \dots &0 & a_1 & b_0
\end{pmatrix}.
\] 
Corresponding orthogonal polynomials $P_n(x)$ are called the persymmetric OPRL.

The persymmetric OPRL can be equivalently characterized by the two equivalent conditions \cite{VZ1}: 

\vspace{4mm}

(i) the weights $w^{(P)}_s$ are given by the formula
\be
w^{(P)}_s = \frac{\sqrt{h_N}}{\left| P_{N+1}'(x_s) \right|} .
\lab{wps} \ee

\vspace{3mm}

(ii) the polynomials $P_{N}(x)$ satisfy the property
\be
P_{N}(x_s) = \sqrt{h_N} (-1)^{N+s} .
\lab{P_N_s} \ee 

Notice that from these conditions one can obtain the following more detailed formula for the persymmetic weights \cite{VZ1}
\be
w_s^{(P)} = (-1)^{N+s} \frac{\sqrt{h_N}}{P_{N+1}'(x_s)}, \quad s=0,1,\dots,N. \lab{w_s_per} \ee

We also mention the following 
\begin{lemma}\label{lem_32}
The spectral points $x_0, x_1, \dots, x_N$ determine the coefficients $b_0, b_1, \dots, b_N$ and $g_1, g_2, \dots g_N$ of the  persymmetric Jacobi matrix $K$ uniquely. 
\end{lemma} 
This lemma follows from formula \re{w_s_per}. Indeed, the weights $w_s^{(P)}$ are determined uniquely form the given spectrum $x_s$. Hence the persymmetric polynomials $P_n^{(P)}(x)$ as well as the entries $g_n, b_n$ of the persymmetic Jacobi matrix are determined uniquely.   This Lemma was first proven in \cite{BG}. It plays an important role in inverse spectral problems for tridiagonal matrices \cite{Gladwell}, \cite{VZ_PST}.

\section{Polynomials orthogonal on the unit circle}
\setcounter{equation}{0}
In this section we recall briefly basic definitions and properties of the polynomials orthogonal on the unit circle (OPUC) introduced by Szeg\H{o} about 100 years ago.

Let $\Phi_n(z)=z^n + O(z^{n-1})$ be a set of monic polynomials with complex coefficients satisfying the recurrence relation \cite{Szego}, \cite{Simon}
\be
\Phi_{n+1}(z) =z \Phi_n(z) - \bar a_n \Phi^*_n(z), \lab{rec_Phi} \ee 
where 
$$
\Phi_n^*(z) =z^n \bar \Phi_n(1/z)
$$
and where $\bar \Phi_n(z)$ being polynomials obtained by complex conjugation of expansion coefficients of the polynomial $\Phi_n(z)$ . The parameters
$$
a_n = -\bar \Phi_{n+1}(0)
$$ 
are important in the theory of OPUC. They are called the Verblunsky parameters \cite{Simon}).

The condition
\be
|a_n|<1, \quad n=0,1,2,\dots , \lab{a_cond} \ee
is equivalent to the statement that the polynomials $\Phi_n(z)$ are orthogonal on the unit circle with respect to a positive measure $d \sigma(\theta)$
\be
\int_0^{2 \pi} \Phi_n(e^{i \theta}) \bar \Phi_m(e^{-i\theta}) d \sigma(\theta) = h_n \delta_{nm}, \lab{ort_Phi} \ee
where
\be
h_n = (1-|a_0|^2)(1-|a_1|^2) \dots (1-|a_{n-1}|^2) \lab{h_n} \ee
are positive normalization constants. 

Similarly to the case of OPRL, one can consider the truncation condition which reduces an infinite family of OPUC to a {\it finite} family of OPUC  $\Phi_1(z), \Phi_2(z), \dots, \Phi_{N+1}(z)$. This truncation condition reads
\be
|a_i|<1, i=0,1,\dots, N-1, \quad a_N= \omega,
\lab{a_omega} \ee
where $\omega = \exp\left(i \chi \right)$ is an arbitrary point on the unit circle: $|\omega|=1$.

Then it is easy to prove (see \cite{Simon}) that the roots of the polynomial $\Phi_{N+1}(z)$ are all distinct and belong to the unit circle:
\be
\Phi_{N+1}(z) = ( z-z_0) (z-z_1) \dots (z-z_N), 
\lab{Phi_N+1} \ee
where $z_k =  \exp\left( i \theta_k \right)$. Moreover, the following condition 
\be
\omega \Phi^*_{N+1}(z) +  \Phi_{N+1}(z) =0
\lab{Phi_Phi} \ee
is valid. 

Define the discrete weights
\be
w_s = \frac{h_N}{\Phi_{N+1}'(z_s) \bar \Phi_N \left( z_s^{-1} \right) }
\lab{w_OPUC} \ee
It is easy to show that these weights are positive $w_s>0, \: s=0,1,\dots, N$ \cite{Simon}. 

Then the following finite orthogonality orthogonality relation \cite{Simon}
\be
\sum_{s=0}^N w_s \Phi_n(z_s) \bar \Phi_m \left( z_s^{-1}\right) = h_n \delta_{nm} 
\lab{fin_ort_OPUC} \ee
is valid for $0 \le n,m \le N$.

It is worth noting that truncation condition \re{a_omega} depends on the choice of the point $\omega$ on the unit circle. By changing $\omega$ we also change the roots $z_s, \: s=0,1,\dots, N$ and the weights $w_s$. Hence depending on the choice of $\omega$ we obtain a family of finite OPUC such that $\Phi_1(z), \dots \Phi_N(z)$ are the same but the final member $\Phi_{N+1}(z; \omega)$ can vary. Such OPUC are called the {\it paraorthogonal} (POPUC) \cite{Gol}.

\section{CMV matrices}
\setcounter{equation}{0}
CMV matrices are unitary operators which play crucial role in spectral theory of the OPUC. Roughly speaking, CMV matrices are analogs of the Jacobi matrices in the case of OPUC.

Given the Verblunsky parameters $a_n$ satisfying conditions \re{a_cond}, one can introduce two semi-infinite block-diagonal matrices 
\be
\mathcal{M}_1 =
 \begin{pmatrix}
1 \\
& a_{1} & \rho_1 &  &    \\
  &\rho_1 & -\bar a_{1} &  &   \\
   & &  &               a_{3} & \rho_3  \\
   & & &              \rho_3 & -\bar a_{3}  \\
  &  & &   & &           a_{5} & \rho_5 \\
   &  & &   & &         \rho_5 & - \bar a_{5}  \\
&   &  & & & & & \ddots  \\
 \end{pmatrix}.
\lab{M1_def} \ee 
and
\be
 \mathcal{M}_2 =
 \begin{pmatrix}
  a_{0} & \rho_0 &  &    \\
  \rho_0 & -\bar a_{0} &  &   \\
   &  &              a_{2} & \rho_2  \\
   &  &              \rho_2 & -\bar a_{2}  \\
  &  &    & &          a_{4} & \rho_4  \\
   &  &    & &         \rho_4 & - \bar a_{4}  \\
&   &  & & & & \ddots  \\
 \end{pmatrix},
\lab{M2_def} \ee
where
\be
\rho_n = \sqrt{1-\left| a_n\right|^2}
\lab{rho_n} \ee
Notice that both matrices are invariant with respect to the reflection over the main diagonal, i.e. they are symmetric  :
\be
{\mathcal{M}_1}^{T} = \mathcal{M}_1, \; {\mathcal{M}_2}^{T} = \mathcal{M}_2
\lab{M1M2_sym} \ee
Moreover, they are both unitary, i.e.
\be
{\mathcal{M}_1}^{\dagger} \mathcal{M}_1 =  {\mathcal{M}_2}^{\dagger}  \mathcal{M}_2 = I
\lab{M1M2_uni} \ee
Hence their product
\be
\mathcal{U}= \mathcal{M}_2 \mathcal{M}_1
\lab{UMM} \ee
is a unitary pentadiagonal "zigzag" operator. This operator plays a crucial role in spectral theory of OPUC (see \cite{Simon} for details). In particular, one can define the following Laurent polynomials \cite{Simon}
\be
\psi_{2n} =\frac{z^{-n} \Phi_{2n}(z)}{\sqrt{h_{2n}}}, \quad \psi_{2n+1} =\frac{z^{n} \bar \Phi_{2n+1}\left(1/z \right)}{\sqrt{h_{2n+1}}}, \: n=0,1,2,\dots
\lab{psi_def} \ee
These polynomials constitute the orthonormal system
\be
\int_0^{2 \pi} \psi_n(e^{i \theta}) \bar \psi_m(e^{-i\theta}) d \sigma(\theta) = \delta_{nm} \lab{ort_psii} \ee
Moreover the column vector
\be
\vec \Psi(z) = \{ \psi_0, \psi_1(z), \psi_2(z), \dots\}^{T}
\lab{vec_psi} \ee
is an eigenvector of the unitary CMV matrix  $\mathcal U$ 
\be
{\mathcal U} \vec \Psi(z) = z \vec \Psi(z).
\lab{eigen_psi} \ee

Truncation condition \re{a_omega} leads to truncation of corresponding CMV matrices. But the "shape" of truncated matrices depends on the parity of $N$.

If $N=2j+1$ is odd then both matrices $\mathcal{M}_{1,2}$ have the {\it even} dimension $N+1=2j+2$ and they have the expressions
\be
\mathcal{M}_1 = {\bf 1} \oplus \Theta_1 \oplus \Theta_3 \oplus \dots \Theta_{2j-1} \oplus {\bf \omega}, \quad 
\mathcal{M}_2 =  \Theta_0 \oplus \Theta_2 \oplus \dots \Theta_{2j} 
\lab{M_tr_e} \ee
where $\Theta_k$ is the $2 \times 2$ matrix
\be
\Theta_k = \begin{pmatrix}
a_k & \rho_k \\
\rho_k  &-\bar a_k
\end{pmatrix}
\lab{theta} \ee
For example, in the simplest case $j=1$ the matrices $\mathcal{M}_1, \mathcal{M}_2$ have the size $4 \times 4$
\be
\mathcal{M}_1 =
 \begin{pmatrix}
1 \\
& a_{1} & \rho_1 &  &    \\
  &\rho_1 & -\bar a_{1} &  &   \\
   & & & \omega
 \end{pmatrix}, \quad \mathcal{M}_2 =
 \begin{pmatrix}
a_{0} & \rho_0 &  &    \\
\rho_0 & -\bar a_{0} &  &   \\
& & a_{2} & \rho_2     \\
& & \rho_2 & -\bar a_{2} .   
 \end{pmatrix}
\lab{M12_4} \ee

If $N=2j$ is even then both matrices $\mathcal{M}_{1,2}$ have the {\it odd} dimension $N+1=2j+1$ and they have the expressions
\be
\mathcal{M}_1 = {\bf 1} \oplus \Theta_1 \oplus \Theta_3 \oplus \dots \Theta_{2j-1} , \quad 
\mathcal{M}_2 =  \Theta_0 \oplus \Theta_2 \oplus \dots \oplus \Theta_{2j-2} \oplus \omega 
\lab{M_tr_0} \ee
In the simplest case $j=1$ both matrices have size $3 \times 3$
\be
\mathcal{M}_1 =
 \begin{pmatrix}
1 & &\\
& a_{1} & \rho_1       \\
&\rho_1 & -\bar a_{1}      
\end{pmatrix}, \quad \mathcal{M}_2 =
 \begin{pmatrix}
a_{0} & \rho_0 &      \\
\rho_0 & -\bar a_{0} &     \\
& & \omega    
 \end{pmatrix}
\lab{M12_3} \ee
The main property of truncated CMV matrices is the following (see \cite{Simon} for details)
\begin{pr}
The spectrum of the unitary matrix $\mathcal{U} = \mathcal{M}_2 \mathcal{M}_1$ coincides with the roots of the polynomial $\Phi_{N+1}(z)$. All these roots are simple and  belong to the unit circle $|z|=1$.
\end{pr}
This proposition is a unitary analog of the corresponding proposition for the OPRL considered in Section 2.  

The eigenvalue problem for truncated CMV matrices reads
\be
{\mathcal U} \vec \Psi_N\left(z_s \right) = z_s \vec \Psi_N\left(z_s \right), \; s=0,1,2,\dots N
\lab{eig_psi_N} \ee
where 
\be
\vec \Psi_N \left(z_s \right) = \left\{\psi_0, \psi_1\left( z_s \right), \dots , \psi_N\left( z_s \right) \right\}^T.
\lab{vec_Psi_N} \ee

\section{Mirror dual OPUC and persymmetric OPUC}
\setcounter{equation}{0}
Let $\Phi_n(z), \: n=0,1,\dots, N$ be a finite set of OPUC such that the Verblunsky parameters satisfy the conditions
\be
a_{-1}=-1, \: |a_n|<1, \: n=0,1,\dots, N-1, \quad a_N=\omega,
\lab{a<1_ini} \ee
where $\omega$ is an arbitrary parameter on the unit circle $|\omega|=1$. Then, as was explained in the previous section, the polynomials $\Phi_n(z)$ satisfy orthogonality condition \re{fin_ort_OPUC} with the weights $w_s$ defined by \re{w_OPUC}.

Define the new set of the Verblunsky parameters \cite{DG1}, \cite{DG2}
 \be
\hat a_n = -\omega \bar a_{N-n-1}.
\lab{hat_a} \ee
It is seen that these parameters satisfy the conditions
\be
\hat a_{-1}=-1, \: |\hat a_n|<1, \: n=0,1,\dots, N-1, \quad \hat a_N=\omega.
\lab{a<1_mirr} \ee
We can define corresponding {\it mirror dual}  OPUC by the recurrence relation
\be
\hat \Phi_0 =1, \; \hat \Phi_{n+1} = z \hat \Phi_n(z) - \hat a_n \hat \Phi_n^*(z)
\lab{rec_hat} \ee
By construction, the mirror dual polynomials $\hat \Phi_n(z)$ satisfy the orthogonality relation
\be
\sum_{s=0}^N \hat w_s \hat \Phi_n\left( \hat z_s \right) {\bar {\hat \Phi}}_m \left( \hat z_s^{-1}\right) = \hat h_n \delta_{nm} 
\lab{mirr_ort_OPUC} \ee
where $\hat z_s$ are the roots of the characteristic polynomial $\hat \Phi_{N+1}(z)$ and where the normalization constants are 
\be
\hat h_0 =1, \; \hat h_n =  \left(1-\left| \hat a_{0}\right|^2 \right) \left(1-\left| \hat a_{1}\right|^2 \right)     \dots \left(1-\left| \hat a_{n-1}\right|^2 \right)
\lab{hat_h} \ee
The weights are defined as 
\be
\hat w_s = \frac{\hat h_N}{\hat \Phi_{N+1}'(\hat z_s) {\bar {\hat \Phi}}_N \left( \hat z_s^{-1} \right) }
\lab{w_mirr} \ee
It is easy to show that \cite{DG1}, \cite{Mar_Si}
\be 
\hat h_{N} = h_N, \quad \hat \Phi_{N+1}(z) = \Phi_{N+1}(z).
\lab{hat_N+1} \ee 
Formula \re{hat_N+1} means that the spectral points $z_0,z_1, \dots, z_{N}$ of the POPUC $\Phi_(z)$ and of $\hat \Phi(z)$ are the same.
Hence the weights $\hat w_s$  can be presented by a more simple formula
\be
\hat w_s = \frac{h_N}{ \Phi_{N+1}'(z_s) {\bar {\hat \Phi}}_N \left(  z_s^{-1} \right) }
\lab{hat_w_1} \ee
In \cite{Mar_Si} it was shown that there is an equivalent presentation of the weights $\hat w_s$:
\be
\hat w_s = \frac{\Phi_N(z_s)}{\Phi_{N+1}'(z_s)}
\lab{hat_w_2} \ee
From these relations we obtain the important 
\begin{pr}
The weights $w_s$ and $\hat w_s$ satisfy the identity
\be
w_s \hat w_s = \frac{h_N}{\left| \Phi_{N+1}'(z_s) \right|^2}.
\lab{ww_id} \ee \end{pr}
Indeed, formula \re{w_OPUC} can be presented in the equivalent form (because $\bar w_s = w_s$)
\be
w_s = \frac{h_N}{\bar \Phi_{N+1}\left( z_s^{-1}\right) \Phi_N(z_s)}.
\lab{w_s_2} \ee
Then multiplying \re{hat_w_2} and \re{w_s_2} we obtain \re{ww_id}. 

Formula \re{ww_id} is an OPUC analog of formula \re{ww_P} for OPRL. 

We call OPUC {\it persymmetric} if the mirror dual sequence $\hat a_n$ coincides with the original sequence $a_n$. In more details this conditions means
\be
a_n =\hat a_n =  -\omega \bar a_{N-1-n}, \: n=0,1,\dots, N.
\lab{per_a} \ee
From \re{per_a} it follows that OPUC are persymmetric iff
\be
\hat \Phi_n(z) = \Phi_n(z), \: n=0,1,\dots, N.
\lab{per_Phi} \ee
Before deriving main properties of the persymmetric OPUC, let us mention restrictions upon the Verblunsky coefficients $a_n$ imposed by the persymmetric conditions \re{per_a}. We fix the parameter $a_N=\omega$ belonging to the unit circle $|\omega|=1$ and  consider conditions for the Verblunsky coefficients $a_0, a_1, \dots, a_{N-1}$.

If $N$ is even then it is sufficient to choose arbitrary $N/2$ complex parameters $a_0, a_1, \dots, a_{(N-2)/2}$ (with the only restriction $|a_k|<1$). Then the rest $N/2$ Varblunsky parameters $a_{N/2}, a_{(N+2)/2}, \dots , a_{N-1}$ are defined uniquely by \re{per_a}. 

If $N$ is odd then again we can choose free $(N-1)/2$ parameters $a_0, a_1, \dots, a_{(N-3)/2}$ (with the same restriction $|a_k|<1$). Then the parameters  $a_{(N+1)/2},a_{(N+3)/2}, \dots a_{N-1}$ are defined uniquely by \re{per_a}. The only exception in this case is the "middle" coefficient $a_{(N-1)/2}$. Indeed, by \re{per_a} it should satisfy the condition 
\be
a_{(N-1)/2} = -\omega \bar a_{(N-1)/2}.
\lab{mid_a_cond} \ee
This condition can be satisfied if one chooses $a_{(N-1)/2}$ in the form
\be
a_{(N-1)/2} = i r \omega^{1/2}, 
\lab{choose_mid_a} \ee
where $r$ is an arbitrary real parameter satisfying the condition $-1 < r <1$. Hence in the case of the odd $N$ we can take $(N-1)/2$ independent complex parameters $a_0, a_1, \dots, a_{(N-2)/2}$ and one free real parameter $r$.

We are ready now to analyze conditions which are equivalent to condition of persymmetric OPUC.

The first such condition is  
\begin{pr}
OPUC are persymmetric iff the weights are
\be 
w_s = \frac{\sqrt{h_N}}{\left| \Phi_{N+1}'(z_s)  \right|}.
\lab{w_per} \ee
\end{pr}
The proof of this proposition follows directly from formula \re{ww_id}: if OPUC are persymmetric then $w_s=\hat w_s$ and formula \re{ww_id} becomes \re{w_per}. Conversely, if the weights $w_s$ satisfy \re{w_per} then from \re{ww_id} it follows that $\hat w_s = w_s$. In turn, this means that $\hat \Phi_n(z) = \Phi_n(z)$.

Another equivalent statement is
\begin{pr}
OPUC are persymmetric iff the values $\Phi_N(z_s)$ lie on the circle with the radius $\sqrt{h_N}$, i.e. iff
\be
|\Phi_N(z_s)| = \sqrt{h_N}, \: s=0,1,\dots N
\lab{Phi_N_sq} \ee
\end{pr}
Again this proposition follows directly from the previous proposition and from formula \re{ww_id}.

Finally, we can present a more detailed description of property \re{Phi_N_sq}. Namely, we have the 
\begin{pr}
Assume that the roots of the polynomial  $\Phi_{N+1}(z)$ of the persymmetric OPUC are
\be
z_s = \exp\left( i \theta_s \right), \: s=0,1,\dots, N,
\lab{z_s_theta} \ee
where $\theta$ are real parameters ordered by increasing:
\be
0 \le \theta_0 < \theta_1 < \dots < \theta_N < 2 \pi
\lab{theta_order} \ee

Then the values of $\Phi_N\left( z_s\right)$ are
\be
\Phi_N\left( z_s\right)= (-1)^{s} \ve \omega^{-1/2} \exp\left( i(N-1) \theta_s/2 \right)\sqrt{h_N} , \: s=0,1, \dots , N, 
\lab{phase_N} \ee
where $\ve = \pm$ is a common factor depending on a choice of an appropriate branch of the value $\omega^{1/2}$.

\end{pr} 

{\it Proof}. Using formula \re{Phi_N_sq} we can present the values of $\Phi_N(z_s)$ as
\be
\Phi_N(z_s) = \left| \Phi_N(z_s) \right| \exp\left(i \chi_s \right) = \sqrt{h_N} \exp\left(i \chi_s \right) ,
\lab{Phi_psi} \ee
where $\chi_s$ is an unknown phase value of $\Phi_N(z_s)$.

By \re{w_OPUC}, we can write
\be
w_s = \frac{h_N}{\Phi_{N+1}'(z_s) \overline{\Phi_n(z_s)}}= \frac{\sqrt{h_N} \exp\left( i \chi_s\right)}{\Phi_{N+1}'(z_s)}.
\lab{w_s_3} \ee
On the other hand, for persymmetric OPUC the weights are given by  \re{w_per}. Comparing this formula with \re{w_s_3} we get
\be
\Phi_{N+1}'(z_s) = \left| \Phi_{N+1}'(z_s)  \right| \exp\left( i \chi_s\right).
\lab{Phi_N+1_exp} \ee
This means that the phase $\exp\left( i \chi_s\right)$ of the complex number $\Phi_N(z_s)$ coincides with the phase of the complex number $\Phi_{N+1}'(z_s)$. It remains to calculate the phase of $\Phi_{N+1}'(z_s)$. 

We have
\be
\Phi_{N+1}'(z_s) = (z_0-z_s) (z_1-z_s) \dots (z_{s-1}-z_s) (z_{s+1}-z_s) \dots (z_N-z_s)
\lab{Phi_N+1_expr} \ee
Substituting to \re{Phi_N+1_expr} the values \re{z_s_theta}, we obtain 
\be
\Phi_{N+1}'(z_s) = \left| \sin\left( \frac{\theta_0-\theta_s}{2} \right) \sin\left( \frac{\theta_1-\theta_s}{2} \right) \dots \sin\left( \frac{\theta_{s-1}-\theta_s}{2} \right) \sin\left( \frac{\theta_{s+1}-\theta_s}{2} \right) \dots \sin\left( \frac{\theta_N-\theta_s}{2} \right) \right| \exp\left( i \chi_n\right),
\lab{N+1_sin} \ee
where 
\be
\exp\left( i\chi_N \right)  = (-1)^{N+s} \exp\left( i \left(\theta_0/2 +\theta_1/2 + \dots + \theta_N/2 \right) \right) \exp\left( i(N-1) \theta_s/2 \right). 
\lab{exp_psi} \ee
On the other hand, 
\be
z_0 z_1  \dots z_N = \exp\left( i \left(\theta_0 +\theta_1 + \dots + \theta_N \right)\right) =(-1)^N \bar a_N = (-1)^N  \omega^{-1}.
\lab{zzz} \ee

This proves our Proposition.

{\it Remark}. The common factor $\ve = \pm 1$ is not essential. Once $\ve$ is chosen then formula \re{phase_N} is well defined.

\vspace{4mm} 

Formula \re{phase_N} can be considered as a OPUC analog of formula \re{w_s_per}. It can be used to create an algorithm to restore uniquely the Verblunsky coefficients of persymmetric polynomials if the spectral points $z_s$ are given.  Indeed, using Lagrange interpolation formula, we can reconstruct the polynomial $\Phi_N(z)$ from its known values at $N+1$ distinct points $z_s$ on the unit circle. Then by inverse Szeg\H{o} formula \cite{Simon} we can reconstruct the polynomials $\Phi_{N-1}(z), \Phi_{N-2}(z), \dots, 1$ uniquely. This automatically generates the Verblunsky coefficients $a_n = -\bar \Phi_{n+1}(0)$. See \cite{per} for corresponding algorithms in the case of OPRL.

\section{Explicit examples}
\setcounter{equation}{0}
Consider first two simplest examples. Let $a_n=0$ for $n=0,1,\dots,N-1$ and $a_N=\omega$ with an arbitrary parameter $\omega$ such that $|\omega|=1$. This is is the so-called "free" case of OPUC \cite{Simon}. It is obvious that
\be 
\Phi_n(z) = z^n, \: n=0,1,\dots, N, \quad \Phi_{N+1}= z^{N+1}-\omega^{-1}
\lab{free_<N} \ee
Putting $\omega = \exp\left(2 i \pi \nu \right)$ with a real parameter $\nu$ we get the spectrum
\be
z_s =\exp\left( \frac{2 i \pi \left(-\nu + s \right)}{N+1}  \right), \; s=0,1,2,\dots , N
\lab{z_s_free} \ee
The weights are equal
\be
w_s = (N+1)^{-1}, \; s=0,1,\dots, N
\lab{w_free} \ee
These OPUC are obviously persymmetric. Property \re{phase_N} can be checked directly. Indeed, for "free" OPUC we have
\be
\Phi_N(z) = z^N.
\lab{Phi_N_free} \ee 
Hence
\be
\Phi_N\left(z_s \right) = z_s^{N} = \frac{z_s^{N+1}}{z_s} = \frac{1}{\omega z_s}
\lab{Phi_N_z_s} \ee
Then formula \re{phase_N} follows easily from \re{Phi_N_z_s}.

The second example corresponds to the so-called single-momentum OPUC \cite{Simon}. They have the Verblunsky parameters 
\be
a_n = -\frac{1}{n+2}, \: n=0,1,2\dots, N-1, \quad a_N=-1. 
\lab{sm_a} \ee
Explicit expression
\be
\Phi_n(z) = \frac{1}{n+1} \: \sum_{k=0}^n (k+1)z^k, \; n=0,1,2,\dots N
\lab{Phi_sm} \ee
and
\be
\Phi_{N+1}(z) = \frac{z^{N+2}-1}{z-1} = z^{N+1} + z^N + \dots + z+1
\lab{sm_N+1} \ee
The normalization coefficients
\be
h_n = \frac{n+2}{2(n+1)} , \: n=0,1,\dots, N
\lab{h_sm} \ee

From \re{sm_N+1} it is seen that the spectrum is
\be
z_s = \exp\left(\frac{2 \pi i (s+1)}{N+2} \right), \; s=0, 1,2,\dots , N
\lab{zs_sm} \ee
I.e. the spectrum includes all roots of unity of order $N+2$ apart from $z=1$.

The corresponding weights are 
\be
w_s  = \left(\frac{N+2}{2} \right) \, \sin^2 \theta_s, \; s=0,1,\dots, N,
\lab{w_sm} \ee
where 
\be
\theta_s = \frac{\pi (s+1)}{N+2}
\lab{theta_sm} \ee
Noe that these weights are normalized
\be
\sum_{s=0}^N w_s =1.
\lab{norm_sm} \ee

What are mirror dual OPUC with respect to the simgle-momentum OPUC? By definition \re{hat_a} we have the Verblunsky parameters
\be
\hat a_n = -\omega a_{N-n-1} = a_{N-n-1} = - \frac{1}{N-n+1}, \: n=0,1,2,\dots N.
\lab{hat_sm_a} \ee
The spectral poits $z_s$ are the same \re{zs_sm}. Remarkably in this case all the weights are equal one to another:
\be
\hat w_s = \frac{1}{N+1}, \: s=0,1,\dots, N
\lab{hat_w_sm} \ee
Formula \re{hat_w_sm} follows directly from \re{ww_id}.

Explicit expression for mirror dual single momentum OPUC is
\be
\hat \Phi_n(z) = z^n + \frac{1}{N-n+2} \left(\frac{z^n-1}{z-1}  \right).
\lab{hat_Phi_sm} \ee
Formula \re{hat_Phi_sm} can be checked directly from recurrence relation \re{rec_Phi}.

What are persymmetric OPUC corresponding to the same spectrum \re{zs_sm}? First of all, we can calculate the weights $w^{(P)}_s$ by using formula \re{w_per}. Indeed, we already know the weights \re{w_sm} for the single momentum OPUC  \re{w_sm}. Formula \re{w_per} then shows that the weights $w^{(P)}_s$ are (up to a normalization factor) square roots of the weights \re{w_sm}. We should be careful only with the normalization factor. Simple calculations lead to the formula
\be
w^{(P)}_s = \tan \nu \, \sin \theta_s
\lab{w_paer_sm} \ee
where
\be
\nu =  \frac{\pi}{2(N+2)}
\lab{nu_def} \ee

with the same $\theta_s$ as in \re{theta_sm}. The weights are appropriately normalized, i.e. 
\be
\sum_{s=0}^N w^{(P)}_s =1.
\lab{norm_wp_sm} \ee

The expression for the Verblunsky parameters is
\be
a_n = -\frac{\sin \nu}{\sin\left(  \nu \left( 2n+3\right) \right)}.
\lab{a_per_sm} \ee
It is seen that $-1< a_n <1, \, n=0,1,\dots, N-1$ and that $a_N=-1$. Formula \re{a_per_sm} follows from results of \cite{Zhe_pol} (Section 4), where several explicit families of polynomials orthogonal on regular polygons were considered.

Finally, we consider the new example of OPUC with the {\it linear} Verblunsky coefficients
\be
a_n = \alpha n + \beta.
\lab{a_lin} \ee
The parameters $\alpha$ and $\beta$ can be found from the following conditions: 

\vspace{4mm}

(i) $a_{-1}=-1$;

\vspace{4mm}

(ii) $a_{N}=\omega$,

where $\omega$ is an arbitrary complex parameters on the unit circle $|\omega|=1$.

The first condition is standard convention in the theory of OPUC \cite{Simon}. The second condition means that we deal with a finite system of OPUC. Indeed, it is obvious that for the linear dependence \re{a_lin} the coefficients $a_n$ cannot satisfy the positivity condition $|a_n|<1$ for all $n$. Hence we need truncation condition like (ii).   From these two conditions we have the expression
\be
a_n = \frac{\omega+1}{N+1} (n+1) -1.
\lab{a_lin_tau} \ee
For all values of the parameter $\omega$ belonging to the unit circle (apart from the degenerate case $\omega=-1$) the Verblunsky coefficients $a_n$ satisfy the standard positivity condition
\be
|a_n|<1, \: n=0,1,\dots, N-1, \; |a_N|=1.
\lab{pos_a_lin} \ee
Indeed, from \re{a_lin_tau} we have
\be
1- \left| a_n \right|^2 = \frac{(\omega+1)^2(n+1)(N-n)}{\omega(N+1)^2}>0, \: n=0,1,\dots, N-1.
\lab{1-Kr} \ee

The case $\omega=-1$ corresponds to the total degeneracy $a_n=-1, \, n=0,1,\dots, N$. Excluding this degenerate case,  we arrive at the case of a finite system of OPUC such that the roots of the polynomials $\Phi_1(z), \Phi_2(z), \dots, \Phi_{N}(z)$ lie inside the unit circle $|z|<1$ while the roots of the polynomial $\Phi_{N+1}(z)$ are all distinct and lie on the unit circle:
\be
\Phi(z_s)=0, \; |z_s|=1, \: s=0,1,\dots, N.
\lab{rt_lin} \ee

Moreover, it is seen from \re{a_lin_tau} that corresponding OPUC are persymmetric:
\be
\hat a_n = \omega \bar a_{N-n-1} = a_n.
\lab{lina_per} \ee

It is possible to present explicit expression of the corresponding OPUC. To do this, let us first introduce the set of the symmetric Krawtchouk polynomials $K_n(x;N)$  defined by the recurrence relation \cite{KLS}
\be
K_{n+1}(x;N) + v_n K_{n-1}(x;N) =xK_n(x;N), \: K_0(x)=1, K_1(x) =x,
 \lab{rec_KrK} \ee
where
\be
v_n = \frac{n(N+2-n)}{4}.
\lab{v_K} \ee
The polynomials $K_{N+2}(x)$ have the roots
\be 
x_k = -(N+1)/2+k, \; k=0,1,\dots, N+1.
\lab{x_kr} \ee 
The Krawtchouk polynomials are orthogonal on the grid of this roots \cite{KLS} 
\be 
\sum_{s=0}^{N+1} K_n(x_s) K_m(x_s) w_s = 0, \quad n \ne m,
\lab{ort_K} \ee
where
\be
w_s = \frac{(N+1)!}{s! (N+1-s)!}
\lab{binomial} \ee
is the binomial distribution.

Define now the orthogonal polynomials $\t K_n(x)$ which are scaled Krawtchouk polynomials, i.e.
\be
\t K_n(x;N) = \kappa^n K_n\left(\frac{x}{\kappa};N \right) ,
\lab{P_K} \ee
where
\be
\kappa^2 =\frac{(\omega+1)^2}{\omega (N+1)}. 
\lab{kappa_tau} \ee
The polynomials $\t K_n(x;N)$ are orthogonal on the scaled grid $\t x_s = \kappa x_s$ with respect to the same binomial weights $w_s$. 

If $\omega = \exp\left(i \sigma \right)$ then all roots $\t x_s$ of the polynomial $\t K_{N+1}(x)$ are located inside the interval $[-2\cos \left(\sigma/2 \right), 2\cos \left(\sigma/2 \right) ]$
\be
-2\cos \left(\sigma/2 \right) \le \t x_0 < \t x_1 < \dots < \t x_{N+1} \le 2\cos\left(\sigma/2 \right).
\lab{int_tau} \ee 
When $\omega=1$ this interval is $[-2,2]$. When $\omega$ is approaching the value $\omega=-1$, the interval is squeezing to zero.  

Define also the  polynomials
\be
P_n(z;N) = z^n \t K_n\left(z+ z^{-1} ;N\right).
\lab{PtK} \ee
These polynomials satisfy the recurrence relation
\be
P_{n+1}(z;N) + \kappa^2 n(N+2-n) z^2 P_{n-1}(z;N) = \left(z^2+ 1  \right) P_n(z;N).
\lab{rec_P_K} \ee
By definition \re{PtK} it is obvious that $P_n(z;N)$ are polynomials of degree $2n$ in the argument $z$. From recurrence relation \re{rec_P_K} it is clear that the polynomials $P_n(z;N)$ are in fact polynomials of degree $n$ in the argument $z^2$.

The polynomials $P_n(z;N)$ are building blocks to express the OPUC explicitly. Namely, we have the 
\begin{pr}
The OPUC with linear Verblunsky coefficients \re{a_lin_tau} have the explicit expression
\be
\Phi_n(z^2;N) = \frac{ P_{n+1}(z;N) -  A_n  P_n(z;N)}{z^2-\omega},
\lab{Phi_P_Kr} \ee
where 
\be
A_n= \frac{\left(\omega +1\right) \left(N -n +1\right)}{N +1}.
\lab{A_K} \ee
\end{pr}
The proof of this proposition can be done by induction. Indeed, it is easy to check that for $n=0$ relation \re{Phi_P_Kr} is valid. Assume that formula \re{Phi_P_Kr} is valid for $n=0,1,\dots, j$. Then for $n=j+1$ formula \re{Phi_P_Kr} follows from recurrence relation \re{rec_P_K}.

Using formula \re{Phi_P_Kr} we can easily find spectral points $z_s$ of the OPUC polynomial $\Phi_{N+1}(z)$ starting from the known spectral points $\t x_s$ of the scaled Krawtchouk polynomials.

\begin{pr}
The roots $z_s =\exp\left( i \theta_s\right)$ of the polynomial $\Phi_{N+1}(z)$ can be found from the transcendent equation
\be
\cos\left( \theta_s/2\right) = \left(\frac{2s}{N+1}-1 \right) \cos\left( \sigma/2\right), \: s=0,1,\dots, N,
\lab{rt_K} \ee
where $\omega =\exp\left( i \sigma\right)$.
The (non-normalized) weights $w_s$ in the orthogonality relation \re{fin_ort_OPUC} are given by the expression 
\be
w_s = \frac{1}{s! (N+1-s)!} \left| \frac{\exp(i\theta_s) -\exp(i \sigma)}{\exp(i\theta_s) -1}\right| = \frac{1}{s! (N+1-s)!} \left|\frac{ \sin\left( \theta_s/2-\sigma/2\right) }{\sin\left(\theta_s/2 \right)} \right| .
\lab{w_K_OPUC} \ee
 \end{pr} 
Notice that if $\omega=1$ then $\sigma=0$ and from \re{rt_K} it is seen that the spectral points cover the whole unit circle: there is the point $z_0=-1$ and the point $z_N=1$. However, when $\sigma \ne 0$ he spectral points avoid the the symmetric arc $|\theta| < |\sigma|$ around the point $z=1$ on the unit circle. When $\omega \to -1$ (i.e. $\sigma \to \pi$) this "forbidden arc" becomes larger and finally covers the whole unit circle, hence in this case all spectral points are located inside of  a small arc near the point $z=-1$.

{\it Proof}. From \re{Phi_P_Kr} it follows that the spectral points $z_s$ (i.e. zeros of $\Phi_{N+1}(z)$ ) correspond to zeros of the polynomial $P_{N+2}\left(y(z) \right)$, where
\be
y(z) = \sqrt{z} + \frac{1}{\sqrt{z}} .
\lab{y_z} \ee
Then formula \re{rt_K} follows. The weights $w_s$ can be derived from the same formula \re{Phi_P_Kr} and from expressions \re{w_per} and \re{wps}  for the weights of the persymmetric polynomials. This leads to formula \re{w_K_OPUC}.

\vspace{4mm}

{\it Remark}. When $\omega=1$ formula \re{Phi_P_Kr} is a special case of the Delsarte-Genin map (DGM) from OPRL to OPUC. This map was introduced in \cite{DG} and developed in \cite{DVZ}, \cite{CMMV}.
Necessary condition for DGM to exist is that the Verblunsky coefficients $a_n$ should be {\it real}. For $\omega \ne 1$ formula \re{Phi_P_Kr} can be considered as a nontrivial generalization of DGM for the case of complex coefficients $a_n$.

\section{Mirror duality relations for CMV matrices}
\setcounter{equation}{0}
In the case of OPRL the mirror dual Jacobi matrix $\hat J$ is connected to the initial Jacobi matrix by the simple relation
\be
RJR = \hat J.
\lab{RJR} \ee
Relation \re{RJR} means that the persymmetric Jacobi matrix $J=\hat J$ is doubly symmetric: with respect to the main diagonal and with respect to the main antidiagonal.   

What is corresponding relation for the CMV matrices? More exactly, we would like to find relation for the  CMV matrix $\mathcal{U}$ similar to \re{RJR} where the matrix $J$ is replaced with the matrix $\mathcal{U}$.

We display here these relations. It happens that, in contrast to the case of OPRL,  relation of type \re{RJR} strongly depends on the parity of $N$ (recall that $N+1$ is the size of quadratic CMV-matrices).

The main tool in our approach will be the so-called {\it quasi-reflection matrices} $Q_N(\tau)$. These matrices are defined as follows. Let $\tau$ be an arbitrary complex parameter belonging to the unit circle, i.e. $|\tau| =1$. Then the matrix $Q_N(\tau)$ is a quadratic  matrix of size $N+1$ which contains only the main antidiagonal with alternating entries $\tau$ and $\bar \tau = \tau^{-1}$:
\be
Q_N(\tau)=\begin{pmatrix}
  0 & 0 & \dots & 0 & \tau    \\
  0 & 0 & \dots  & 
 \tau^{-1} & 0  \\
    \dots  & \dots & \dots & \dots & \dots      \\
   \tau^{\pm 1} & 0 &  \dots & 0 &0  \\
\end{pmatrix}.
\lab{Q_def} \ee 
The value $\tau$ or $\tau^{-1}$ on the left bottom corner of $Q_N(\tau)$ depends on parity of $N$. If $N$ is odd, the size of the matrix $Q_N(\tau)$ is even and this "final" value is $\tau^{-1}$. If $N$ is even then the final value is $\tau$. For example, for $N=3$ we have
\be
Q_3(\tau)=\begin{pmatrix}
  0 & 0 & 0 &\tau    \\
  0 & 0 & \tau^{-1}  & 0 \\
 0 & \tau & 0 & 0 \\
 \tau^{-1} & 0 & 0 & 0
\end{pmatrix}.
\lab{Q_3} \ee  
Similarly, for $N=4$ the quasi-reflection matrix is
\be
Q_4(\tau)=\begin{pmatrix}
  0 & 0 & 0 & 0 & \tau    \\
  0 & 0 & 0& \tau^{-1}  & 0 \\
 0 & 0 & \tau & 0 & 0 \\
 0& \tau^{-1} & 0 & 0 & 0\\
 \tau & 0 & 0 & 0 & 0
\end{pmatrix}.
\lab{Q_4} \ee  
The following two propositions describe main properties of the quasi-reflection matrices.

\begin{pr}
If $N$ is odd then the quasi-reflection matrix is an involution, i.e. $Q_N^2(\tau)={I}$; the eigenvalues of the matrix $Q_N(\tau)$ are $\pm 1$. The number $(N+1)/2$ of the eigenvalues $1$ and $-1$ is the same.
\end{pr}

For the case of even $N$ the operators $Q_N(\tau)$ have different properties. Namely, we have the
\begin{pr}
If $N$ is even then the  the quasi-reflection matrix is not an involution, i.e.  $Q_N^2(\tau) \ne {I}$ (until $\tau \ne \pm 1$). Nevertheless, the matrix $Q_N(\tau)$ is unitary, i.e. $Q_N(\tau) Q_N^{\dagger}(\tau)=I$ (as usual, the symbol $Q^{\dagger}$ means the Hermitian conjugate matrix). Moreover, $Q_N^{\dagger}(\tau) = Q_N\left( \tau^{-1}\right)$, so that $Q_N(\tau)Q_N\left( \tau^{-1}\right)=I $. The spectrum of the matrix $Q_N(\tau)$ consists of 4 values: $\pm \tau$ and $\pm \tau^{-1}$.
\end{pr}

The proof of these Propositions is elementary and we skip them. 

We can now apply the quasi-reflection matrices to describe invariance properties of the mirror dual CMV matrices.

Consider first a more simple case of {\it odd} vlue of $N$.

We have the
\begin{pr}
Let $N$ is odd and $\mathcal{M}_1$ and $\mathcal{M}_2$ be CMV matrices of even size $N+1$ defined by \re{M_tr_e}. Introduce the mirror dual CMV matrices $\mathcal{ \hat M}_1$ and $\mathcal{\hat M}_2$ with entries $\hat a_n$ defined by \re{hat_a}. Then the following relations take place
\be
\mathcal{ \hat M}_1 = Q_{N}\left(\tau^{-1} \right) \mathcal{M}_1 Q_N \left( \tau \right), \quad  \mathcal{ \hat M}_2 = Q_{N}\left(\tau \right) \mathcal{M}_2 Q_{N}\left( \tau^{-1} \right),
\lab{hat_LM} \ee
where $\tau=\omega^{-1/2}$.
\end{pr} 
Proof of this proposition is direct.

Consider an example $N=5$. Then 
\be
\mathcal{M}_1 =
 \begin{pmatrix}
1 & & & & &\\
& \bar a_{1} & \rho_1 & &       \\
&\rho_1 & - a_{1} & &  \\
& &  & \bar a_{3} & \rho_3 &   \\
& &  & \rho_3 & -a_3  &  \\
& & & & & \omega^{-1}
\end{pmatrix}, \quad \mathcal{M}_2 =
 \begin{pmatrix}
\bar a_{0} & \rho_0 &  & &     \\
\rho_0 & - a_{0} &  & &    \\
& & \bar a_2 & \rho_2 & \\
& &   \rho_2 & -a_2 & \\
& & & & \bar a_4 & \rho_4 \\
& & &  &  \rho_4 & -a_4
 \end{pmatrix}.
\lab{M12_5} \ee
Mirror reflected CMV matrices are
\be
\mathcal{\hat M}_1 =
 \begin{pmatrix}
1 & & & & &\\
& -\omega^{-1} a_{3} & \rho_3 & &       \\
&\rho_3 & \omega \bar a_{3} & &  \\
& &  & -\omega^{-1} a_{1} & \rho_1 &   \\
& &  & \rho_1 &  \omega \bar a_1  &  \\
& & & & & \omega^{-1}
\end{pmatrix}, \quad \mathcal{\hat M}_2 =
 \begin{pmatrix}
-\omega^{-1}  a_{4} & \rho_4 &  & &     \\
\rho_4 & \omega \bar a_{4} &  & &    \\
& & -\omega^{-1} a_2 & \rho_2 & \\
& &   \rho_2 & \omega \bar a_2 & \\
& & & & -\omega^{-1} a_0 & \rho_0 \\
& & &  &  \rho_0 &  \omega \bar a_0
 \end{pmatrix}.
\lab{M12_5M} \ee
Then it is easy to check relations \re{hat_LM} directly.

Now we can construct the unitary operator ${\mathcal U} ={\mathcal M}_2 {\mathcal M}_1$. For the mirror reflected unitary operator we have
\be
\hat U =  \mathcal {\hat M}_2 \mathcal {\hat M}_1 = Q_{N}\left(\tau \right) \mathcal{M}_2 Q_N \left( \tau^{-1} \right) Q_{N}\left(\tau^{-1} \right) \mathcal{M}_1 Q_{N}\left( \tau \right) = Q_{N}\left( \tau \right) \mathcal{U}  Q_{N}\left( \tau \right).
\lab{hat_U_odd_N} \ee
Indeed, the operator $Q_N \left( \tau^{-1} \right)$ is an involution: $Q_N \left( \tau^{-1} \right)Q_N \left( \tau^{-1} \right)=I$. This yields \re{hat_U_odd_N}.

Hence we have the
\begin{pr}
For $N$ odd the CMV unitary operator is transformed as
\be
 \mathcal {\hat U} = Q_{N}\left( \tau \right) \mathcal{U}  Q_{N}\left( \tau \right), \quad \mbox{where} \quad  \tau = \omega^{-1/2}.
\lab{UU_odd_N} \ee
\end{pr}
We can now compare relation \re{UU_odd_N} with corresponding relation \re{RJR} for Jacobi matrices. It is seen that formally both relations look similarly. The only distinction is that the involution operator $Q_N(\tau)$ has slightly more complicated structure \re{Q_def}. In particular, it is seen that \re{UU_odd_N} is similarity transformation and hence the spectra of both matrices $\mathcal U$ and $\mathcal {\hat U}$ are the same.

Consider the case when $N$ is even. In this case we have the 
\begin{pr}
 For $N$ even the following relations take place
\be
\mathcal{ \hat M}_2 = Q_{N}(\tau^{-1} ) \mathcal{M}_1 Q_N \left( \tau^{-1} \right), \quad  \mathcal{ \hat M}_1 = Q_{N}\left(\tau \right) \mathcal{M}_2 Q_{N}\left( \tau \right),
\lab{hat_LM_N_even} \ee
where $\tau=\omega^{1/2}$. This means that the unitary CMV operator is transformed as
\be
Q_N(\tau) {\mathcal U} Q_N \left( \tau^{-1} \right) = {\mathcal M}_1 {\mathcal M}_2 = \mathcal{ \hat U}^{T}.
\lab{tr_U_N_even} \ee
\end{pr} 
Formula \re{tr_U_N_even} shows that the operator $\mathcal U$ is unitary equivalent to the operator $\mathcal{ \hat U}^{T}$. This again means that the spectra of the operators $\mathcal U$ and $\mathcal {\hat U}$ are the same. However, the main distinction with respect to the case of $N$ odd is that the operator $\mathcal U$ is transformed to the transposed operator $\mathcal{ \hat U}^{T}$.

Consider now action of the quasi-reflection operators $Q_N(\tau)$ on the vectors $\vec v_s = \vec \Psi_N(z_s), \: s=0,1,2,\dots ,N$ defined by \re{vec_Psi_N}, where $z_s= \exp\left(i \theta_s \right)$ are roots of the polynomial $\Phi_{N+1}(z)$ ordered by increasing of $\theta_s$ (see \re{theta_order}).

We already know that the vectors $\vec v_s, : s=0,1,\dots, N$ are the eigenvectors of the unitary CMV operator
\be
{\mathcal U} \vec v_s = z_s \vec v_s, \; s=0,1,\dots,N
\lab{UVs} \ee
with distinct eigenvalues $z_s$. The action of the operator $Q(\tau)$ on the vectors $\vec v_s$ is given by the simple
\begin{pr}
The vectors 
\be
{\hat {\vec v}}_s = Q_N(\tau) \vec v_s, \; s=0,1,\dots, N
\lab{hat_v} \ee
are eigenvectors of the unitary CMV operator $\hat {\mathcal U}$ with the same eigenvalues:
\be
\hat {\mathcal U} {\hat {\vec v}}_s = z_s {\hat {\vec v}}_s, \; s=0,1,\dots, N
\lab{hat_UV} \ee
\end{pr}
Proof of this proposition is elementary. Indeed, 
\be
\hat {\mathcal U} {\hat {\vec v}}_s = Q_N(\tau) \mathcal U Q_N(\tau) Q_N(\tau) \vec v_s = Q_N(\tau) U \vec v_s = z_s {\hat {\vec v}}_s,
\lab{hat_UV_proof} \ee
where we have used the involution property of the operator $Q_N(\tau)$.

Consider now the special case when the system of OPUC is persymmetric and $N$ is still odd. This means that $\hat {\mathcal U} = \mathcal U$. The operator $Q_N(\tau)$ commutes with $\mathcal U$ in this case:
\be
{\mathcal U} Q_N(\tau) =  Q_N(\tau) {\mathcal U}.
\lab{comm_QU} \ee
In turn, this means that the vector ${\hat {\vec v}}_s$ should coincide with the vector $\vec v_s$ to within a factor:
\be
{\hat {\vec v}}_s = \mu_s \vec v_s, \; s=0,1,\dots, N
\lab{hat_vv} \ee
Equivalently, this means that the vector $\vec v_s$ is an eigenvector of the operator $Q_N(\tau)$ with the eiegnevalue $\mu_s$
\be
Q_N(\tau) \vec v_s = \mu_s \vec v_s.
\lab{Q_vv} \ee
On the other hand, we already know that all eigenvalues of the operator $Q_N(\tau)$ are $ \mu_s=\pm 1$. In more details,we have the 
\begin{pr}
If $N$ is odd and the spectrum of the operator $\mathcal U$ is ordered by increasing of the aguemnts $\theta_s$ \re{theta_order}, then 
\be
Q_N(\tau) \vec v_s = \ve (-1)^s \vec v_s, \: s=0,1,\dots, N,
\lab{Q_pm} \ee
where $\ve$ is a common factor $\pm 1$ like in \re{phase_N}.
\end{pr}
The proof of this proposition follows from Proposition 5.4 and from definition \re{psi_def} of the Laurent polynomials $\psi_n(z)$.

From this proposition we have an important 

{\bf Corollary}

{\it For $N$ odd there is the relation}
\be 
 \psi_{N-n}\left(z_s \right) = \ve (-1)^s \omega^{\nu_n} \psi_{n}\left(z_s \right), \; n,s=0,1,\dots, N
\lab{psi_psi_rel} \ee
{\it where $\nu_n=1/2$ for $n$ even and $\nu_n =-1/2$ for $n$ odd}.

This Corollary follows directly from the action of the antidiagonal matrix $Q_N(\tau)$ on the components $\psi_n\left(z_s\right)$ of the vector $\vec v_s$. The case $n=0$ is equivalent to formula \re{phase_N}.

\vspace{4mm}

{\it Remark}. There is similar property for the persymmetric  OPRL \cite{Bor}, \cite{per}. However, for the persymmetric OPUC there are two main distinctions: 1) formula \re{psi_psi_rel} is valid only for $N$ odd; 2) there are factors $\omega^{\pm 1/2}$ in rhs of \re{psi_psi_rel} which are absent in the case of OPRL.

\section{Conclusion}
\setcounter{equation}{0}
We have demonstrated that the persymmetric OPUC have many properties similar to the persymmetric OPRL. For example, the discrete weight function of the persymmetric OPUC is determined by the spectrum of CMV matrix uniquely. This means that the inverse spectral problem for CMV matrices has  unique solution if the spectrum is given. This result is similar to the corresponding results in the theory of persymmetric Jacobi matrices.

However, CMV matrices have nontrivial features which have no analogs in the theory of the persymmetric Jacobi matrices. Namely, when the size of the persymmetric CMV matrix is even, then it commutes with a reflection operator. So in this cse there is strong similarity with the case of persymmetric Jacobi matrices (note that the reflection operator is different than the trivial reflection operator in the case of persymmetric Jacobi matrix) . When the size of CMV matrix is odd, then such matrix does not commute with (quasi)reflection matrices.

This property means that properties of the persymmetric CMV matrices are more interesting and nontrivial. The most intriguing question is: what are possible applications  of the persymmetric CMV matrices? It is well known that persymmetric Jacobi matrices play a crucial role in the theory of the perfect state transfer \cite{Albanese}, \cite{VZ_PST}. One can expect that persymmetric CMV matrices could play similar role in transferring the quantum information as well. This is an interesting open problem.  

Another wide area of future research is constructing explicit examples of the persymmetric OPUC. We have presented in this paper two nontrivial examples of the persymmetric OPUC: the first one is a persymmetric analog of the single-moment OPUC and the second one is a circle analog of the Krawtchouk polynomials (in this case the OPUC contain an additional complex parameter which is absent for the Krawtchouk OPRL). For OPRL we know many explicit examples of persymmetric polynomials and corresponding Jacobi matrices. In particular, any "classical" grid - linear, quadratic and hyperbolic - generates explicit persymmetric OPRL \cite{VZ_PST}. We expect that for the unit circle there are appropriate "classical" grids which generate corresponding persymmetric OPUC.

\bb{99}

\bi{AG} N.I. Akhiezer and I.M.Glazman, {\it Theory of linear operators in Hilbert space}, 2nd rev.ed., "Nauka",
Moscow, 1966;   English transl. of 3rd ed., Pitman, Boston, MA, 1981. Reprint by Dover, 1993.

\bi{Albanese} C.~Albanese, M.~Christandl, N.~Datta and A.~Ekert, {\it Mirror inversion of quantum states in linear registers}, Phys. Rev. Lett. {\bf 93} (2004), 230502.

\bi{Atkinson} F.~V.~Atkinson, {\it Discrete and Continuous Boundary problems}, Academic Press, NY, London, 1964.

\bi{BG} C.~de~Boor and G.~H.~Golub {\it The numerically stable reconstruction of a Jacobi matrix from spectral data}, Lin. Alg. Appl. {\bf 21} (1978), 245--260.

\bi{BS} C.~de~Boor and E.~Saff, {\it Finite sequences of orthogonal polynomials connected by a Jacobi matrix}, Lin. Alg. Appl. {\bf 75} (1986), 43--55.

\bi{Bor} A.~Borodin, {\it Duality of Orthogonal Polynomials on a
Finite Set}, J. Stat. Phys. {\bf 109}, (2002), 1109--1120.

\bi{CMV} M.J. Cantero, L. Moral, and L. Vel\'azquez, {\it Five-diagonal matrices and zeros
of orthogonal polynomials on the unit circle}, Linear Algebra Appl. {\bf 362} (2003),
29--56.

\bi{CMMV} M. J. Cantero, F. Marcell\'an, L. Moral, L. Vel\'azquez, {\it A CMV connection between orthogonal polynomials on
the unit circle and the real line}, J.Approx.Theory, {\bf 266}, (2021), 105579.

\bi{Chi} T.~Chihara, {\it An Introduction to Orthogonal
Polynomials}, Gordon and Breach, NY, 1978.




\bi{DG} P. Delsarte, Y. Genin, {\it The split Levinson algorithm}, IEEE Trans. Acoust.
Speech Signal Process. {\bf 34} (1986), 470--478.

\bi{DG1} P. Delsarte, Y. Genin,  {\it Tridiagonal approach to the algebraic environment of
Toeplitz matrices. I. Basic results}, SIAM J. Matrix Anal. Appl. {\bf 12} (1991), no. 2,
220--238.

\bi{DG2} P. Delsarte, Y. Genin, {\it Tridiagonal approach to the algebraic environment of
Toeplitz matrices. II. Zero and eigenvalue problems}, SIAM J. Matrix Anal. Appl.
{\bf 12} (1991), no. 3, 432--448.

\bi{DVZ} M.Derevyagin, L.Vinet, A.Zhedanov, {\it CMV matrices and Little and Big -1
Jacobi Polynomials}, Constr. Approx. {\bf 36} (2012), 513--535.

\bi{per} V.Genest, S.Tsujimoto, L.Vinet and A.Zhedanov, {\it Persymmetric Jacobi matrices, isospectral deformations and orthogonal polynomials}, J.Math.anal.Appl., {\bf 450}, 915--928 (2017). arXiv:1605.00708. 

\bi{Gladwell}  G.~M.~L.~Gladwell , {\it Inverse Problems in Vibration}, 400 pp., Martinus Noordhoff, 1986.

\bi{Gol} L.Golinskii, {\it Quadrature formula and zeros of para-orthogonal polynomials on the unit circle}, Acta Mathematica Hungarica {\bf 96}, 169--186 (2002).

\bi{GoKu} L.Golinskii and M.Kudryavtsev, {\it An inverse spectral theory for finite CMV matrices},   
Inverse Problems and Imaging {\bf 4}, 93--110 (2010).  arXiv:0705.4353 (2007).

\bi{GoKu2} L.Golinskii and M.Kudryavtsev, {\it Rational interpolation and mixed inverse spectral problem for finite CMV matrices}, J.Approx.Theory, {\bf 159}, 61--84, (2009).


\bi{Ismail} M. E. H. Ismail, {\it Classical and Quantum Orthogonal Polynomails in One Variable}, Encyclopedia of
Mathematics and its Applications {\bf 98}, Cambridge University Press, 2009.

\bi{KLS} R.~Koekoek, P.~Lesky and R.~Swarttouw, {\it Hypergeometric Orthogonal Polynomials and Their Q-analogues}, Springer-Verlag, 2010.


\bi{Mar_Si} A. Mart\'inez-Finkelshtein, B. Simanek, and B. Simon. {\it Poncelet's theorem, paraorthogonal polynomials and the
numerical range of compressed multiplication operators}. Adv. Math., {\bf 349}, 992--1035, (2019).


\bi{Simon} B.Simon, {\it Orthogonal Polynomials On The Unit
Circle}, AMS, 2005.

\bi{Szego} G.~Szeg\H{o}, 
{\it Orthogonal Polynomials}, fourth edition,  AMS, 1975.

\bi{VZ1} L.~Vinet and A.~Zhedanov, 
{\it A characterization of classical and semiclassical orthogonal polynomials from their dual polynomials},  
J. Comp. Appl. Math. {\bf 172} (2004), 41--48.

\bi{VZ_PST} L.~Vinet and A.~Zhedanov, 
{\it How to construct spin chains with perfect state transfer}, 
Phys. Rev. A {\bf 85} (2012), 012323.

\bi{Zhe_pol} A.Zhedanov, {\it On the Polynomials Orthogonal on Regular Polygons}, J.Approx. Theory {\bf 97}, 1--14 (1999).

\eb

\end{document}